\documentclass{amsart}

\usepackage[table,xcdraw]{xcolor}

\usepackage{amsfonts, amsmath, amssymb, amsthm, blkarray, color, enumerate, float, graphicx, longtable, hyperref, mathtools, nicematrix, pgfplots, soul, stmaryrd, tikz-cd, todonotes, verbatim} 

\hypersetup{colorlinks=true, linkcolor=blue,  anchorcolor=blue, citecolor=blue, filecolor=blue, menucolor=blue, pagecolor=blue, urlcolor=blue}

\usepackage[all]{xy}
\usepackage[outline]{contour}

\usepackage{tikz-cd}
\usepackage{pgfplots}

\usetikzlibrary{arrows, calc, fit, intersections, matrix,decorations.pathreplacing}

\usepackage[normalem]{ulem}
\newtagform{black}[\color{black}]{(}{)}

\usepackage{nicematrix}

\tikzset{
    ncbar angle/.initial=90,
    ncbar/.style={
        to path=(\tikztostart)
        -- ($(\tikztostart)!#1!\pgfkeysvalueof{/tikz/ncbar angle}:(\tikztotarget)$)
        -- ($(\tikztotarget)!($(\tikztostart)!#1!\pgfkeysvalueof{/tikz/ncbar angle}:(\tikztotarget)$)!\pgfkeysvalueof{/tikz/ncbar angle}:(\tikztostart)$)
        -- (\tikztotarget)
    },
    ncbar/.default=0.5cm,
}

\newcounter{markeq}
\newcommand{\pstrut}[1]{\vrule height0pt depth0pt width0pt #1 \fboxsep}
\newcommand*\bmarkeq{\stepcounter{markeq}%
  \tikz[remember picture]\node(startframe-\themarkeq){\pstrut{height}};%
  \kern\fboxsep}
\newcommand*\emarkeq{\kern\fboxsep
  \begin{tikzpicture}[remember picture,overlay]
    \node (endframe-\themarkeq){\pstrut{depth}};
    \draw[,red,opacity=0.8] (startframe-\themarkeq.north) 
      rectangle (endframe-\themarkeq.south);
  \end{tikzpicture}%
}

\tikzset{square left brace/.style={ncbar=0.5cm}}
\tikzset{square right brace/.style={ncbar=-0.5cm}}

\tikzset{round left paren/.style={ncbar=0.5cm,out=120,in=-120}}
\tikzset{round right paren/.style={ncbar=0.5cm,out=60,in=-60}}

\newcount\colveccount
\newcommand*\colvec[1]{
        \global\colveccount#1
        \begin{pmatrix}
        \colvecnext
}
\def\colvecnext#1{
        #1
        \global\advance\colveccount-1
        \ifnum\colveccount>0
                \\
                \expandafter\colvecnext
        \else
                \end{pmatrix}
        \fi
}

\makeatletter
\renewcommand*\env@matrix[1][*\c@MaxMatrixCols c]{%
  \hskip -\arraycolsep
  \let\@ifnextchar\new@ifnextchar
  \array{#1}}
\makeatother
\newcommand\restr[2]{{
  \left.\kern-\nulldelimiterspace 
  #1 
  \vphantom{\big|} 
  \right|_{#2} 
  }}

\allowdisplaybreaks
\makeatletter
\newcommand{\leqnomode}{\tagsleft@true}
\newcommand{\reqnomode}{\tagsleft@false}
\makeatother

\theoremstyle{definition}

\newtheorem{remark}{Remark}

\allowdisplaybreaks

\makeatletter
\@addtoreset{footnote}{page}
\makeatother

\author[C. Glass]{Cheyne Glass}
\address{Cheyne Glass, State University of New York at New Paltz, Department of Mathematics, 1 Hawk Dr., New Paltz, NY 12561}
  \email{glassc@newpaltz.edu}

\author[E. Vidaurre]{Elizabeth Vidaurre}
\address{Elizabeth Vidaurre, Molloy University,  1000 Hempstead Avenue, Rockville Centre, NY 11570}
  \email{evidaurre@molloy.edu}

\title{Topological Data Analysis\\ via Undergraduate Linear Algebra}

\keywords{Linear Algebra, Topological Data Analysis, Persistent Homology, Undergraduate Mathematics}

\begin{document}

\begin{abstract}
Topological Data Analysis has grown in popularity in recent years as a way to apply tools from algebraic topology to large data sets. One of the main tools in topological data analysis is persistent homology. This paper uses undergraduate linear algebra to provide explicit methods for, and examples of, computing persistent (co)homology. 
\end{abstract}
\maketitle

\tableofcontents

\section{A very brief introduction to TDA}

Given a data set in the form of a point cloud, perhaps the scatterplot for a spreadsheet of numerical values, we're often interested in describing the perceived geometry of this point cloud. If the dimension of our point cloud (i.e. perhaps the number of columns of the spreadsheet) is at most three, then we can simply look at a scatter plot of the data points in $\mathbb{R}^n$ to visualize our point cloud and try to qualitatively assign geometric features. For instance, does it form clusters? Does it have a void or hole? These questions are not easily answered by the traditional statistical methods in an undergraduate curriculum. Using topology to study point clouds gives us new tools. Topology is an area of mathematics which is uniquely qualified to study the general shape of a point cloud without being too concerned with details such as angles, lengths, or linearity. In this paper, we explicitly share one of the topological tools developed to study data: \emph{persistent homology}. It has  proven itself useful in various applications, including the study of data in materials science of glass \cite{glass} and in neuroscience \cite{Be}, among others.

To illustrate how topology can be helpful, consider some examples of $2$-dimensional point clouds in Figure \ref{Fig:RegressionPlots} below.

\begin{figure}[H]
   \begin{minipage}{0.45\textwidth}
     \centering
     \includegraphics[width=\textwidth]{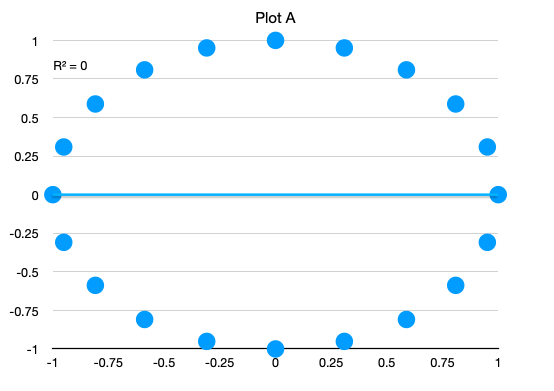}
   \end{minipage}
   \begin{minipage}{0.45\textwidth}
     \centering
     \includegraphics[width=\textwidth]{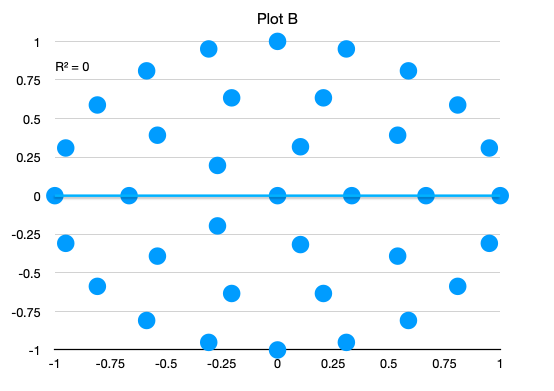}
   \end{minipage}
      \begin{minipage}{0.45\textwidth}
     \centering
     \includegraphics[width=\textwidth]{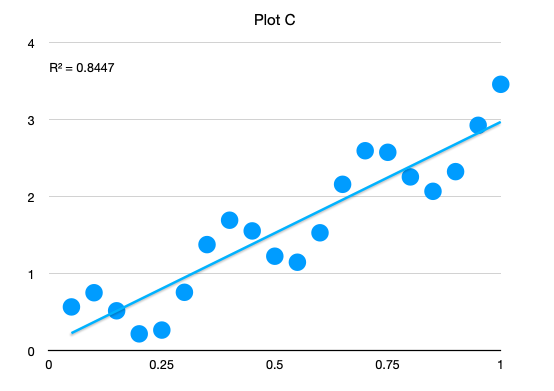}
   \end{minipage}
        \caption{Scatterplots $A$, $B$, and $C$ with $R^2$ values of $0$, $0$, and $0.8447$, respectively}\label{Fig:RegressionPlots}
\end{figure}
Through the lens of any usual form of regression (one of the main tools usually taught in undergraduate statistics or linear algebra), there would be no difference between plots $A$ and $B$, despite the fact that a quick glance would leave no doubt that these two data sets have very different structures. Studying these plots' {persistent homology} \cite{CZ} (or persistence modules) fails to distinguish plots $B$ and $C$. We mention this failure, where standard regression analysis  succeeds, in order to make it clear that topological data analysis is not a replacement of regression techniques but an additional tool.

We offer a friendly introduction to persistent homology through \v{Cech} cohomology. In Sections \ref{PC} and \ref{advanced}, we discuss how to extract topological information from a data set without requiring prior knowledge of topology concepts. Section \ref{computation} details the main ``by hand'' computational example.

\section{The cohomology of a space}\label{PC}
In this section, we will  show how to construct a space from a data set and then the cohomology of that space.

In general, for a point cloud of data, $\mathcal{D} \subset \mathbb{R}^n$, we wish to take this discrete set of data points, with no interesting topology, and associate a smooth/continuous space with an interesting topology, for which we can think of our data points as sampling from. We begin by choosing a certain resolution value $\epsilon > 0$, and treat the points in $\mathcal{D}$ as Euclidean balls of radius $\epsilon$. The union of these balls will give us some associated topological space we can refer to as $\mathcal{D}_{\epsilon} \subset \mathbb{R}^n$ (see Figure \ref{fig: Simple Point Cloud} for a visual). For a fixed $\epsilon$, we can then look at the topological invariants of this space. In this paper, we will focus on the \v{C}ech Cohomology\footnote{The choice of homology or cohomology, and which construction one uses to compute the (co)homology is largely a matter of taste at the early stages of one's learning about topological data analysis. The authors found it quicker and more student friendly to explain this construction as opposed to taking the homology of the related simplicial complex. However, it should be noted that in practice cohomology is more efficient: \url{https://mathoverflow.net/questions/290226/why-is-persistent-cohomology-so-much-faster-than-persistent-homology}} of the spaces $\mathcal{D}_{\epsilon}$.

The \v{Cech} cohomology of $\mathcal{D}_{\epsilon}$ is a sequence of vector spaces, $\check{H}_{\epsilon}^0, \check{H}_{\epsilon}^1, \check{H}_{\epsilon}^2, \ldots$ which will be able to distinguish the connected components, holes, etc, of $\mathcal{D}_{\epsilon}$ just as discussed above for the plot in Figure \ref{Fig:RegressionPlots}. For example, in the case of Plot $A$ above (for certain choices of $\epsilon$ which ``connects the dots'') the fact that it seems to be a single connected cluster will correspond to $\check{H}_{\epsilon}^0$ being one-dimensional while the fact that it has a single non-trivial loop will correspond to $\check{H}_{\epsilon}^1$ being one dimensional. Plots $B$ and $C$, on the other hand, are perceived to be connected but have no non-trivial loops, and as a result, the dimensions of their corresponding vector spaces have values  $dim\left(\check{H}_{\epsilon}^0\right) = 1$ and  $dim\left(\check{H}_{\epsilon}^1\right)= 0$. The dimensions of these cohomology vector spaces are often referred to in algebraic topology as the \emph{Betti}\footnote{The validity of this sentence actually depends on a beautiful symmetry within the Betti numbers coming from \emph{Poincare duality}, which we will happily not worry about trying to explain.} numbers, denoted $\beta_i = dim(\check{H}_{\epsilon}^i)$, of our topological spaces $\mathcal{D}_{\epsilon}$. In theoretical algebraic topology, one is quite happy with the Betti numbers (i.e. dimensions) alone since they do not depend on some parameter like $\epsilon$, but for the purposes of \emph{persistence} we will need to keep track of how these vector spaces vary under different choices of epsilon.

\subsection{The cohomology for a choice of $\epsilon$ via example}

Consider the following data set, 
\begin{table}[h]
\caption{Simple Data Set with a ``hole''}
\label{my-label}
\begin{tabular}{lllllllllllllllllllll}
\rowcolor[HTML]{C0C0C0} 
Point  & x-value & y-value \\
\cellcolor[HTML]{C0C0C0}$p_0$ & 1 & 2\\
\cellcolor[HTML]{C0C0C0}$p_1$ & 2 & 3\\
\cellcolor[HTML]{C0C0C0}$p_2$ & 2 & 1\\
\cellcolor[HTML]{C0C0C0}$p_3$ & 3.7 & 2  \\
\cellcolor[HTML]{C0C0C0}$p_4$ & 3.7 & 4.5  \\
\end{tabular}
\end{table}
and plot each point on an $xy$-plane. We will study the set (or the space) given by placing a ball of radius $\epsilon$ around each point:

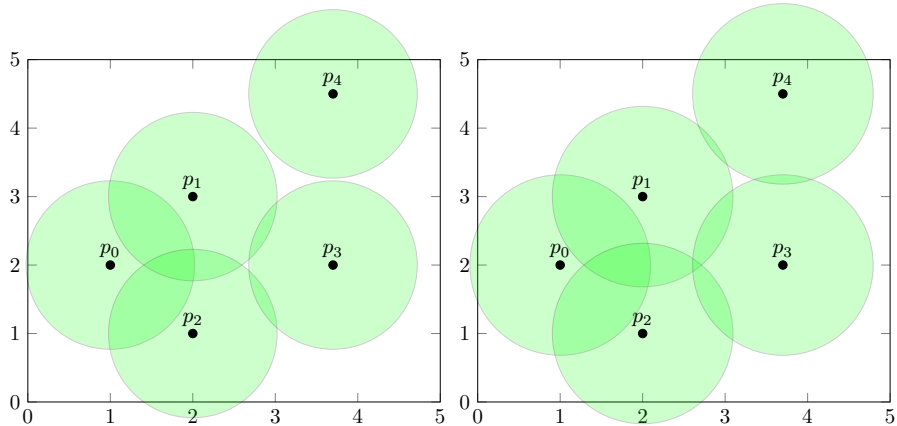
\begin{figure}[H]
   \begin{center}
\begin{tikzpicture}[scale=0.8]
\begin{axis}[xmin=0, xmax=5, ymin=0, ymax=5,nodes near coords={$p_\coordindex$}]
\addplot[only marks,mark=*,mark size=1.4cm,fill=green,opacity=0.2] coordinates {(1,2)  (2,3)  (2,1) (3.7, 2) (3.7, 4.5)};
\addplot[only marks,mark=*,mark size=2pt,fill=black]               coordinates {(1,2)  (2,3)  (2,1) (3.7, 2) (3.7,4.5)};
\end{axis}
\end{tikzpicture}\begin{tikzpicture}[scale=0.8]
\begin{axis}[xmin=0, xmax=5, ymin=0, ymax=5,nodes near coords={$p_\coordindex$}]
\addplot[only marks,mark=*,mark size=1.5cm,fill=green,opacity=0.2] coordinates {(1,2)  (2,3)  (2,1) (3.7, 2) (3.7, 4.5)};
\addplot[only marks,mark=*,mark size=2pt,fill=black]               coordinates {(1,2)  (2,3)  (2,1) (3.7, 2) (3.7,4.5)};
\end{axis}
\end{tikzpicture}
\caption{The collection of points is the set $\mathcal{D} = \{p_0, p_1, p_2, p_3, p_4\}$, and the radius of each disk centered at the points is our $\epsilon$. The figure on the left uses $\epsilon = 1.4$ and the figure on the right uses $\epsilon' = 1.5$. The union of the disks is the topological space, $\mathcal{D}_{\epsilon}$ so that on the left we have $\mathcal{D}_{\epsilon}$ and on the right we have $\mathcal{D}_{\epsilon'}$.}\label{fig: Simple Point Cloud}
\end{center}
\end{figure}

Focusing on the choice above of $\epsilon =1.4$, we can see that there is a triple intersection between the $\epsilon$-balls centered around points $p_0, p_1$, and $p_2$, but the three $\epsilon$-balls centered around the points $p_1, p_2$, and $p_3$ are only pairwise intersecting. We can also see that the $\epsilon$-ball around $p_4$ does not intersect any others. We will now show how to obtain the two \emph{Betti numbers}, $\beta_0$ and $\beta_1$, for this data-set-turned-topological-space. 

\subsection{The vector spaces in the complex}
To begin, we set up a sequence of vector spaces and linear maps between them based upon our data set, $\mathcal{D}$, and our choice of $\epsilon$, 
\begin{equation}\label{EQ:CechComplex}
 0 \xrightarrow{\delta^{-1}} \check{C}^{0}_{\epsilon}(\mathcal{D})\xrightarrow{\delta^{0}} \check{C}^{1}_{\epsilon}(\mathcal{D})\xrightarrow{\delta^{1}} \check{C}^{2}_{\epsilon}(\mathcal{D}) \to \dots
 \end{equation}
which we will refer to as the \emph{\v{C}ech Complex} \cite[Section 10]{BT}, \cite[Section 1.3]{Br}. Each vector space $\check{C}^{p}_{\epsilon}(\mathcal{D})$ is determined by the $(p+1)$-fold intersections of the $\epsilon$-balls centered at the points of $\mathcal{D}$. To be explicit, 
\begin{itemize}
\item $\check{C}^{0}_{\epsilon}(\mathcal{D})$ is the (free) vector space whose basis corresponds to the labelings of the five $\epsilon$-balls:
\[ \{u_0, u_1, u_2, u_3, u_5\}, \] and so is a $5$-dimensional vector space isomorphic to $\mathbb{R}^5$.
\item $\check{C}^{1}_{\epsilon}(\mathcal{D})$ is the (free) vector space whose basis corresponds to the labelings of the five different double-intersections of those four $\epsilon$-balls:
\[ \{u_{01}, u_{02}, u_{12}, u_{13}, u_{23}\}, \] and so is a $5$-dimensional vector space isomorphic to $\mathbb{R}^5$. Notice, for instance, that $u_{01}$ is included because the $\epsilon-$balls corresponding around the points 0 and 1 are overlapping, whereas there is no $u_{03}$. 
\item $\check{C}^{2}_{\epsilon}(\mathcal{D})$ is the (free) vector space whose basis corresponds to the labelings of the triple intersection(s):
\[ \{u_{012}\}, \] and so is a $1$-dimensional vector space isomorphic to $\mathbb{R}^1$.
\item $\check{C}^{p}_{\epsilon}(\mathcal{D})= \{0\}$, for $p>2$ is the trivial vector space as there are no higher intersections. Similarly, we define all of the vector spaces, $\check{C}^{p}_{\epsilon}(\mathcal{D})$,  for $p<0$ to be trivial.
\end{itemize}

\subsection{The linear maps in the complex}
Next, we describe the linear maps between the vector spaces $ \check{C}^{p}_{\epsilon}(\mathcal{D})$, called \emph{boundary maps}, for our \v{C}ech complex:
\[\check{C}^{p}_{\epsilon}(\mathcal{D}) \xrightarrow{\delta^{p}} \check{C}^{p+1}_{\epsilon}(\mathcal{D}).\]
Again, we will try to unpack the formal definition provided for example in \cite[Section 1.3]{Br}. The components of vectors in $\check{C}^{p}_{\epsilon}(\mathcal{D})$ are essentially (locally constant) choices of real numbers assigned to $p$-fold intersections of open sets. Thus, a linear map like $\delta^p$ above must assign a real number to each $(p+1)$-fold intersection of open sets, using the numbers assigned to all of the $p$-fold intersections involved. So, given a vector
\[ f_{\bullet} = \begin{pmatrix} f_0 \\ f_1 \\ f_2 \\ f_3 \\ f_4\end{pmatrix} \in \check{C}^{0}_{\epsilon}(\mathcal{D})\]
whose components $f_i$ are numerical values chosen for the corresponding epsilon ball centered at $p_i$, we see that its image must be of the form
\[ \delta^0\left(f_{\bullet} \right)= \begin{pmatrix} g_{01} \\ g_{02} \\ g_{12} \\ g_{13}\\ g_{23} \end{pmatrix} \in \check{C}^{1}_{\epsilon}(\mathcal{D}),\]
where each component $g_{ij}$ of its image under $\delta^0$, corresponds to the intersection of the epsilon balls centered at $p_i$ and $p_j$. The way we are defining these $g_{ij}$ components in this complex is merely by writing $g_{ij} = f_j - f_i$. This definition matches with the formal one given in \cite{Br}. 

For our particular \v{C}ech complex, built from our four epsilon balls and their intersections, working out the above assignments results in the following matrix:
\begin{equation}\label{EQ:delta0example}
\delta^0 = \NiceMatrixOptions{code-for-last-col = \color{blue},code-for-first-row = \color{blue}} \begin{bNiceMatrix}[first-row, last-col]
u_0 & u_1 & u_2 & u_3 & \\
-1&  1 & 0 & 0 & u_{01}\\
-1&  0 & 1 & 0 & u_{02}\\
0&  -1 & 1 & 0 & u_{12}\\
0&  -1 & 0 & 1 & u_{13}\\
0&  0 & -1 & 1 & u_{23}\\
\end{bNiceMatrix}
\end{equation}

For the second boundary map, $\check{C}^{1} \xrightarrow{\delta^{1}} \check{C}^{2}$, we proceed in an analogous fashion to the above. Given a vector in $\check{C}^{1}$, whose components $g_{ij}$ are corresponding to the $2$-fold intersections, each component $h_{ijk}:=g_{jk} - g_{ik} + g_{ij}$ (again here we are unpacking the referenced general definition for the reader) of the image of $\delta^1$ in $\check{C}^{2}$ corresponds to the $3$-fold intersections. Thus we obtain the matrix, 
\begin{equation}\label{EQ:delta1example}
\delta^1 = \NiceMatrixOptions{code-for-last-col = \color{blue},code-for-first-row = \color{blue}}
\begin{bNiceMatrix}[first-row, last-col]
u_{01} & u_{02} & u_{12} & u_{13} & u_{23} &  \\
1&  -1 & 1 & 0 &  0 & u_{012}\\
0&  0 & 0 & 0 &  0 & u_{123}
\end{bNiceMatrix}
\end{equation}

At this point, we have the following \v{C}ech complex:
\[ \dots \to 0 \xrightarrow{\delta^{-1}_{\epsilon}=0} \check{C}^{0}_{\epsilon}\cong \mathbb{R}^5 \xrightarrow{\delta^{0}_{\epsilon}} \check{C}^{1}_{\epsilon} \cong \mathbb{R}^5 \xrightarrow{\delta^{1}_{\epsilon}} \check{C}^{2}_{\epsilon} \cong \mathbb{R}^1 \xrightarrow{\delta^2 = 0} \check{C}^3_{\epsilon}\simeq \mathbb{R}^0 \to \dots \]

\subsection{The Betti numbers of the complex}

Finally, we compute the \emph{Betti numbers} by taking the difference of the dimensions of the kernel and image of the boundary maps at each stage (see page \pageref{PAGE:cohomology} for more details):
\[ \beta^p := dim \left( ker(\delta^p)\right) -  dim \left( im(\delta^{p-1})\right)\]

In order to compute the dimensions we must first compute the rank, $dim \left( im(\delta^p)\right)$), of the boundary maps (by counting the pivots after applying your favorite row-reduction algorithm, for example) and then applying the Rank-Nullity theorem. For our boundary maps above we have:
\[ dim \left( im(\delta^{0})\right) =3 \quad \text{and} \quad dim \left( im(\delta^{1})\right) = 1\]
with all other (trivial) boundary maps having rank equal to zero. By the Rank-Nullity theorem, we can then find the dimensions of the kernels:
\begin{align*}
dim \left( ker(\delta^0)\right) &=dim \left( \check{C}^{0} \right) - dim \left( im(\delta^0)\right)= 5-3 = 2.\\
dim \left( ker(\delta^1)\right) &=dim \left( \check{C}^{1} \right) - dim \left( im(\delta^1)\right)= 5-1 = 4.\\
dim \left( ker(\delta^2)\right) &=dim \left( \check{C}^{2} \right) - dim \left( im(\delta^2)\right)= 1-0 = 1.
\end{align*}

Finally, we can compute our \emph{Betti} numbers:
\begin{align*}
\beta^0 &= dim \left( ker(\delta^0)\right) -  dim \left( im(\delta^{-1})\right) = 2- 0 = 2.\\
\beta^1 &= dim \left( ker(\delta^1)\right) -  dim \left( im(\delta^{0})\right) = 4- 3 = 1.\\
\beta^2 &= dim \left( ker(\delta^2)\right) -  dim \left( im(\delta^{1})\right) = 1- 1 = 0.
\end{align*}

A topologist would look at these three Betti numbers, and make three very different inferences. The $0$-th Betti number, $\beta^0 = 2$ tells the topologist that this space is not \emph{connected}, but has two connected components: the one component corresponding to the cluster around $p_0, p_1, p_2$, and $p_3$ and the other component corresponding to the solitary cluster around $p_4$. In terms of a data cloud, this is suggesting that they are \emph{clustered} around two groups of data points. The $1$st Betti number, $\beta^1 = 1$ is loosely saying that there is at least one non-trivial loop in the space. This would tell the topologist (and the data analyst) that the space (of data points) is forming some sort of \emph{hole}. In our particular case, this is precisely due to that very small gap between data points $p_1, p_2$, and $p_3$ in Figure \ref{fig: Simple Point Cloud}, which didn't allow a triple intersection to form there. While the appearance of a ``hole'' in this example might seem an exaggeration, the difference between $\beta_1$ being equal to 1 or 0 is precisely the difference between Plot A versus Plots B and C in Figure \ref{Fig:RegressionPlots}. Finally, the $2$nd Betti number, $\beta_2 = 0$, was mostly due to the fact that our data was two-dimensional, even if we had formed higher-intersections. If it were three dimensional, and $\beta_2 \ne 0$, this would suggest a spherical hole present in our data\footnote{An object with a spherical hole would be like a basketball as opposed to a billiard ball}.

For more inspiration as to what the Betti numbers would tell a data analyst, we encourage the reader to browse their favorite references\footnote{Don't think too hard: \url{https://en.wikipedia.org/wiki/Betti_number}} to see what the Betti numbers of higher dimensional spaces, like tori or spheres, are.


%


In this section, we computed the cohomology of $D_\epsilon$ for a fixed $\epsilon$. However, to study a big data set, we want to let $\epsilon$ vary so that we can see which features, such as loops, really stand out or ``persist''. While cohomology is a very important tool in topology, \emph{persistent (co)homology} is one of the topological tools used to study data sets. The final section of this paper will consider this issue of persistence.

\subsection{Existing computational tools for implementation} 

To do these computations with a real data set, which would usually have more than five points, there are various ways to do them by computer. There is a browser based program, \href{https://live.ripser.org}{live.ripser.org}, that outputs the \emph{barcode} (or persistence diagram). A barcode is a visual representation of the persistent homology showing a set of intervals (or thin bars) above the real line for each dimension. Each line represents a topological feature, starts at the value of $\epsilon$ where that feature appeared, and ends at the value of $\epsilon$ at which that feature no longer exists. For instance, in dimension one, the bar code for Plot A in Figure  \ref{Fig:RegressionPlots} would have one long interval, whereas a point cloud that resembled a figure eight with one circle larger than the other would have two intervals, one longer than the other.

For those with more coding experience, there are also  existing computer programs that compute persistent homology. For tools in  python, see \href{https://scikit-tda.org/}{scikit-tda.org/}, and in R, \href{https://cran.r-project.org/web/packages/TDA/index.html}{cran.r-project.org/web/packages/TDA/index.html}.

\section{Advanced considerations: measuring persistence} \label{advanced}
This section collects some concepts from linear algebra, including quotient spaces, which are necessary to understand what is happening when we study persistence. We will end with some explicit calculations in our final section that are small enough to do by hand (as the authors and their students have!). 

\subsection{Betti numbers redux}\label{PAGE:cohomology}
Given a \v{C}ech complex 
\[  0 \xrightarrow{\delta^{-1}} \check{C}^{0}_{\epsilon}(\mathcal{D})\xrightarrow{\delta^{0}} \check{C}^{1}_{\epsilon}(\mathcal{D})\xrightarrow{\delta^{1}} \check{C}^{2}_{\epsilon}(\mathcal{D}) \to \dots \]
as in equation \eqref{EQ:CechComplex}, what we really would like to compute are its \emph{(co)homology vector spaces}\footnote{More generally these are usually only carrying the structures of \emph{graded modules} or \emph{groups}, but in our case at each level they are vector spaces.}, denoted here as $\check{H}^p_{\epsilon} \left( \mathcal{D} \right)$. These 
vector spaces 
are defined \cite[Section 6]{W} by taking the quotient space of the kernel of the outgoing map by the image of the incoming map:
\[ \check{H}^p_{\epsilon} \left( \mathcal{D} \right):= ker \left( \delta^p \right) /  im \left( \delta^{p-1} \right).\]
For the purposes of the Betti numbers for a single choice of $\epsilon$, computing the dimension of this quotient vector space is enough and is given by  the difference of the dimensions of $ker \left( \delta^p \right)$ and $im\left( \delta^{p-1} \right)$. However, in practice, keeping track of the bases which generate these quotient spaces is quite relevant to understanding which topological features \emph{persist} as we vary the resolution $\epsilon$. Thus, understanding this quotient more carefully by a coherent choice of bases is related to the computation of a data set's persistent homology.

\subsection{Letting the resolution vary}\label{PAGE:resolution_varies}
Given a single data set, $\mathcal{D}$, but two choices of resolution, $\epsilon  < \epsilon'$, we can study how the cohomology of the \v{C}ech complex transforms. Recall that our \v{C}ech complexes were constructed by considering $p$-many points' $\epsilon$-disks all intersecting. But notice that if the disks of the smaller radius, $\epsilon$, centered at $p$-many data points are all intersecting, then we know the disks of the larger radius, $\epsilon'$, are also forming a $p$-fold intersection. This fairly intuitive phenomenon produces a \emph{chain map}, which is a sequence of linear maps $f^{p}: \check{C}^{p}_{\epsilon'} \to \check{C}^{p}_{\epsilon}$ satisfying the condition $f^p \cdot \delta^p_{\epsilon'} = \delta^p_{\epsilon} \cdot f^{p-1}$ for each $p = 0, 1, 2, \ldots$. A common way to visualize a chain map is to think of it as a sequence of vertical arrows, 
\begin{equation}\label{EQ:chain map diagram}
\begin{tikzcd}
0 \arrow{r} & \check{C}^{0}_{\epsilon'} \arrow{r}{\delta^0_{\epsilon'}} \arrow{d}{f^0} &  \check{C}^{1}_{\epsilon'} \arrow{d}{f^1} \arrow{r}{\delta^{1}_{\epsilon'}} &  \dots  \arrow{r}{\delta^{p-1}_{\epsilon'}}  & \check{C}^{p}_{\epsilon'} \arrow{r}{\delta^p_{\epsilon'}} \arrow{d}{f^{p-1}} &  \check{C}^{p+1}_{\epsilon'} \arrow{d}{f^p} \arrow{r}{\delta^{p+1}_{\epsilon'}} & \dots \\
0 \arrow{r} &  \check{C}^{0}_{\epsilon} \arrow{r}{\delta^0_{\epsilon}} &  \check{C}^{1}_{\epsilon}  \arrow{r}{\delta^{1}_{\epsilon}} & \dots  \arrow{r}{\delta^{p-1}_{\epsilon}}  & \check{C}^{p}_{\epsilon} \arrow{r}{\delta^p_{\epsilon}} &  \check{C}^{p+1}_{\epsilon}  \arrow{r}{\delta^{p+1}_{\epsilon}} & \dots 
\end{tikzcd}
\end{equation}
whereby traveling through each path of compositions gives the same result (referred to as a \emph{commutative diagram}). This chain map in turn produces a linear map between vector spaces representing the cohomology, $\check{H}^p \left( \check{C}^{\bullet}_{\epsilon'}\right) \to \check{H}^p \left( \check{C}^{\bullet}_{\epsilon}\right)$,
for each $p$ \cite{W}. Features like connected components, and holes, are represented by basis vectors for the quotient spaces. These basis vectors are special  \emph{cohomology classes} in these vector spaces.  If a class is sent to zero under one of the above linear maps for some change in resolution, then that class is said to be ``(co-)born'' at this new resolution. If a class is not in the image of such linear map then it is said to ``(co-)die'' at this new resolution. Tracking the ``birth and death'' of these topological features is a key aspect of persistent (co)homology \cite{C}, as it allows us to see how long a certain features ``lives'' or persists. The longer a feature persists, the more likely it is to have some real-world meaning. 

In order to fully approach this persistence mathematically then, one would need to algorithmically compute the following:
\begin{enumerate}[(i)]
\item Construct the \v{C}ech complex for a maximal sampling of resolutions, $\epsilon$.
\item Compute the cohomology (quotient) vector spaces for each dimension $p$ and each resolution $\epsilon$ by a clever use of the bases, so that for each change in resolution, they can track which classes are born and which die. 
\item Present the persistence of the classes in a digestible manner. 
\end{enumerate}

Thankfully, others (see \cite{CZ}, for example, for a collection of such references) have done incredible work to program all of the above. But that doesn't mean there isn't a lot to learn by trying to do it by hand as we now show below.

\section{An explicit computation of persistence}\label{computation}
In this final section, we offer an explicit example (computed by hand) of the ideas discussed in this paper. It should be noted that the example used in this paper was highly curated for simplicity, and the authors have generally used slightly more involved examples in talks or assignments, so the reader should feel free to believe that they can get pretty far by hand or with some modest computational tools.

Recall that for our choice of $\epsilon =1.4$ for  $\mathcal{D}_{\epsilon=1.4}$ from Figure \ref{fig: Simple Point Cloud}, we constructed the complex, 
\[ 0 \xrightarrow{\delta^{-1}_{\epsilon}=0} \check{C}^{0}_{\epsilon}\cong \mathbb{R}^5 \xrightarrow{\delta^{0}_{\epsilon}} \check{C}^{1}_{\epsilon} \cong \mathbb{R}^5 \xrightarrow{\delta^{1}_{\epsilon}} \check{C}^{2}_{\epsilon} \cong \mathbb{R}^1 \xrightarrow{\delta^2_{\epsilon} = 0} \check{C}^3_{\epsilon} \simeq \mathbb{R}^0 \to \dots \]
with bases:
\begin{itemize}
\item $\check{C}^{0}_{\epsilon} = \text{span}\left\{ u_0, u_1, u_2, u_3, u_4 \right\}$,
\item $\check{C}^{1}_{\epsilon} = \text{span}\left\{ u_{01}, u_{02}, u_{12}, u_{13}, u_{23} \right\}$, and
\item $\check{C}^{2}_{\epsilon} = \text{span}\left\{ u_{012} \right\}$.
\end{itemize} 
and nontrivial matrices, $\delta^0_{\epsilon}$ from \eqref{EQ:delta0example} and $\delta^1_{\epsilon}$ from \eqref{EQ:delta1example}. Similarly, by going through the same considerations but for $\epsilon' = 1.5$, we can construct the complex, 
\[ 0 \xrightarrow{\delta^{-1}_{\epsilon'}=0} \check{C}^{0}_{\epsilon'}\cong \mathbb{R}^5 \xrightarrow{\delta^{0}_{\epsilon'}} \check{C}^{1}_{\epsilon'} \cong \mathbb{R}^7 \xrightarrow{\delta^{1}_{\epsilon'}} \check{C}^{2}_{\epsilon'} \cong \mathbb{R}^2 \xrightarrow{\delta^2_{\epsilon'} = 0} \check{C}^3_{\epsilon'} \dots \]
with bases:
\begin{itemize}
\item $\check{C}^{0}_{\epsilon'} = \text{span}\left\{ u'_0, u'_1, u'_2, u'_3, u'_4 \right\}$,
\item $\check{C}^{1}_{\epsilon'} = \text{span}\left\{ u'_{01}, u'_{02}, u'_{12}, u'_{13},  u'_{14}, u'_{23}, u'_{34} \right\}$,  and
\item $\check{C}^{2}_{\epsilon'} = \text{span}\left\{ u'_{012}, u'_{123} \right\}$.
\end{itemize} 
The chain maps $f^p$ are defined in a fairly intuitive way which we will describe implicitly by:
\begin{itemize}
\item $f^0: \check{C}^{0}_{\epsilon'} \to \check{C}^{0}_{\epsilon}$ sends each basis vector $u'$ to its corresponding basis vector $u$. 
\item $f^1: \check{C}^{1}_{\epsilon'} \to \check{C}^{1}_{\epsilon}$ sends each basis vector $u'$ whose index has a ``$4$'' to the zero vector and all other basis vectors to their corresponding basis vector $u$.
\item  $f^2: \check{C}^{2}_{\epsilon'} \to \check{C}^{2}_{\epsilon}$ sends $u'_{012}$ to $u_{012}$ and sends $u_{123}$ to zero\footnote{The inquisitive reader might notice the maps $f^p$ could also be defined naturally in the opposite direction which is true, but for ``\v{C}ech reasons'' the maps we use fit into a larger discussion.}. 
\end{itemize}
It is important to check that the induced linear maps satisfy the condition stated above \eqref{EQ:chain map diagram};  this a fun exercise  with equations of matrices. 
\subsection{A basis for the cohomology induced by $\epsilon'$}
The only two non-trivial boundary maps $\delta^p_{\epsilon'}$ are now shown below in \eqref{EQ:delta0primereduction} and \eqref{EQ:delta1primereduction}  as an (vertically-)augmented matrix, along with a specific method of reduction which hopefully will become clearer as we proceed. The easier of the calculations is to compute $\check{H}^0_{\epsilon'} =  ker\left( \delta^0_{\epsilon'} \right) / im\left( \delta^{-1}_{\epsilon'} \right)$, which is isomorphic to  $ker\left( \delta^0_{\epsilon'} \right)$ since  $im\left( \delta^{-1}_{\epsilon'} \right)\simeq\mathbb{R}^0$, by column reduction:

\begin{align}
\begin{bNiceMatrix}
\delta^0_{\epsilon'}\\
\hline
Id_{\check{C}^0_{\epsilon'}}
\end{bNiceMatrix} = &\NiceMatrixOptions{code-for-first-row = \color{blue},
code-for-last-col = \color{blue}}
\begin{bNiceMatrix}[first-row, last-col]
u'_0 & u'_1 & u'_2 & u'_3 & u'_4 & \\
-1&  1 & 0 & 0  & 0 & u'_{01}\\
-1&  0 & 1 & 0  & 0 & u'_{02}\\
0&  -1 & 1 & 0  & 0 & u'_{12}\\
0&  -1 & 0 & 1  & 0 & u'_{13}\\
0&  -1 & 0 & 0  & 1 & u'_{14}\\
0&  0 & -1 & 1  & 0 & u'_{23}\\
0&  0 & 0 & -1  & 1 & u'_{34}\\
\hline
1&  0 & 0 & 0  & 0 & u'_{0}\\
0&  1 & 0 & 0  & 0 & u'_{1}\\
0&  0 & 1 & 0  & 0 & u'_{2}\\
0&  0 & 0 & 1  & 0 & u'_{3}\\
0&  0 & 0 & 0  & 1 & u'_{4}
\end{bNiceMatrix} \xmapsto{\substack{col_2 +=  col_1}}
\begin{bNiceMatrix}
-1&  0 & 0 & 0  & 0 \\
-1&  -1 & 1 & 0  & 0 \\
0&  -1 & 1 & 0  & 0 \\
0&  -1 & 0 & 1  & 0 \\
0&  -1 & 0 & 0  & 1 \\
0&  0 & -1 & 1  & 0 \\
0&  0 & 0 & -1  & 1\\
\hline
1&  1 & 0 & 0  & 0 \\
0&  1 & 0 & 0  & 0\\
0&  0 & 1 & 0  & 0\\
0&  0 & 0 & 1  & 0 \\
0&  0 & 0 & 0  & 1 
\end{bNiceMatrix}\label{EQ:delta0primereduction}\\
\xmapsto{\substack{col_3 +=  col_2}}
&\begin{bNiceMatrix}
-1&  0 & 0 & 0  & 0 \\
-1&  -1 & 0 & 0  & 0 \\
0&  -1 & 0 & 0  & 0 \\
0&  -1 & -1 & 1  & 0 \\
0&  -1 & -1 & 0  & 1 \\
0&  0 & -1 & 1  & 0 \\
0&  0 & 0 & -1  & 1\\
\hline
1&  1 & 1 & 0  & 0 \\
0&  1 & 1 & 0  & 0\\
0&  0 & 1 & 0  & 0\\
0&  0 & 0 & 1  & 0 \\
0&  0 & 0 & 0  & 1 
\end{bNiceMatrix} \xmapsto{\substack{col_4 +=  col_3}}
\begin{bNiceMatrix}
-1&  0 & 0 & 0  & 0 \\
-1&  -1 & 0 & 0  & 0 \\
0&  -1 & 0 & 0  & 0 \\
0&  -1 & -1 & 0  & 0 \\
0&  -1 & -1 & -1  & 1 \\
0&  0 & -1 & 0  & 0 \\
0&  0 & 0 & -1  & 1\\
\hline
1&  1 & 1 & 1  & 0 \\
0&  1 & 1 & 1  & 0\\
0&  0 & 1 & 1  & 0\\
0&  0 & 0 & 1  & 0 \\
0&  0 & 0 & 0  & 1 
\end{bNiceMatrix}\\
\xmapsto{\substack{col_5 +=  col_4}}
&\begin{bNiceMatrix}[last-col]
-1&  0 & 0 & 0  & 0 & u'_{01}\\
-1&  -1 & 0 & 0  & 0& u'_{02}\\
0&  -1 & 0 & 0  & 0& u'_{12}\\
0&  -1 & -1 & 0  & 0& u'_{13}\\
0&  -1 & -1 & -1  & 0& u'_{14}\\
0&  0 & -1 & 0  & 0  &u'_{23}\\
0&  0 & 0 & -1  & 0& u'_{34}\\
\hline
1&  1 & 1 & 1  & 1 & u'_{0}\\
0&  1 & 1 & 1  & 1& u'_{1}\\
0&  0 & 1 & 1  & 1& u'_{2}\\
0&  0 & 0 & 1  & 1 & u'_{3}\\
0&  0 & 0 & 0  & 1 & u'_{4}
\CodeAfter
\tikz \draw[rounded corners] (8-5) -| (last-|6) -- (last-|5) |- (8-5) ;
\end{bNiceMatrix}=\begin{bNiceMatrix}
\delta^0_{\epsilon'} \cdot L\\
\hline
Id_{\check{C}^0_{\epsilon'}} \cdot L = L
\end{bNiceMatrix} \label{EQ: H0 epsilon prime reduction} 
\end{align}

The single column of zeros on the top block of the reduced matrix point us to the circled column below it as the single basis vector, $u'_0 + u'_1 + u'_2 + u'_3 + u'_4$, for the null space yielding:
\begin{align}\label{EQ:BasisKerDelta0prime}
\check{H}^0_{\epsilon'} \simeq \ker\left( \delta^0_{\epsilon'} \right) &= \text{span} \left\{u'_0 + u'_1 + u'_2 + u'_3 + u'_4 \right\} . 
\end{align}

Next, we want to reduce the matrix for $\delta^1_{\epsilon'}$ but we must take care: if we are to compute the quotient spaces $\check{H}^p_{\epsilon'} = ker\left( \delta^p_{\epsilon'} \right)  / im\left( \delta^{p-1}_{\epsilon'} \right)$ efficiently we want to ensure that the basis for the column space includes into the basis for the null space. Before being explicit we outline the procedure \cite{K} abstractly:
\begin{remark}\label{REM: KerIm bases together}
Given two matrices $A$ and $B$ satisfying $0 = AB$, we note that  $im(B) \subset ker(A)$. Performing column reductions on $A$ to yield $A \cdot Q$, we can perform induced\footnote{Switching between these reduction matrices and their inverses we find to be a very nice synthesis of reduction concepts learned in class.} row reductions via $Q^{-1} \cdot B$, preserving $0 = (A \cdot Q) \cdot (Q^{-1} \cdot B)$. However, $Q^{-1} \cdot B$ will only be ``partially reduced''. So we further row-reduce it via $P \cdot Q^{-1} \cdot B$, and then now perform the induced column reductions via $A \cdot Q \cdot P^{-1}$, again preserving $0 = (A \cdot Q \cdot P^{-1} ) \cdot (P \cdot Q^{-1} \cdot B)$. This last step will not alter the kernel subspace, but will ensure that the basis we get for the kernel is compatible with the image basis we obtain for $B$; i.e. the image basis will now be a subset of this kernel basis.  
\end{remark}
Below are some highlights from the computations. While we personally performed every row/column operation by hand, we ommitted those individual steps merely for brevity in this paper\footnote{See Cheyne Glass' \href{https://cheynejglass.com/2022/11/05/an-elementary-algorithm-for-persistent-cohomology/}{blog post} for the full details.}: 

\begin{align}\label{EQ:delta1primereduction}
\NiceMatrixOptions{code-for-first-row = \color{blue},
code-for-last-col = \color{blue}}
\begin{bNiceMatrix}
\delta^1_{\epsilon'}\\
\hline
Id_{\check{C}^1_{\epsilon'}}
\end{bNiceMatrix} =& 
\NiceMatrixOptions{code-for-first-row = \color{blue},
code-for-last-col = \color{blue}}\begin{bNiceMatrix}[first-row, last-col]
u'_{01} & u'_{02} & u'_{12} & u'_{13} & u'_{14} & u'_{23} & u'_{34} &  \\
1&  -1 & 1 & 0 & 0 & 0  & 0 & u'_{012}\\
0&  0 & 1 & -1 & 0 & 1  & 0 & u'_{123}\\
\hline
1&  0 & 0 & 0 & 0 & 0  & 0 & u'_{01}\\
0&  1 & 0 & 0 & 0 & 0  & 0 &u'_{02}\\
0&  0 & 1 & 0 & 0 & 0  & 0 &u'_{12}\\
0&  0 & 0 & 1 & 0 & 0  & 0 &u'_{13}\\
0&  0 & 0 & 0 & 1 & 0  & 0 &u'_{14}\\
0&  0 & 0 & 0 & 0 & 1  & 0 &u'_{23}\\
0&  0 & 0 & 0 & 0 & 0  & 1 &u'_{34}\\
\end{bNiceMatrix}  \leadsto \begin{bNiceMatrix}
1&  0 & 0 & 0 & 0 & 0  & 0 \\
0&  0 & 1 & 0 & 0 & 0  & 0 \\
\hline
1&  1 & -1 & -1 & 0 & 1  & 0 \\
0&  1 & 0 & 0 & 0 & 0  & 0 \\
0&  0 & 1 & 1 & 0 & -1  & 0 \\
0&  0 & 0 & 1 & 0 & 0  & 0 \\
0&  0 & 0 & 0 & 1 & 0  & 0 \\
0&  0 & 0 & 0 & 0 & 1  & 0 \\
0&  0 & 0 & 0 & 0 & 0  & 1 \\
\end{bNiceMatrix} = \begin{bNiceMatrix}
\delta^1_{\epsilon'} \cdot Q\\
\hline
Q
\end{bNiceMatrix}
\end{align}
Solving for $Q^{-1}$ by either inverting the column operations or some other preferred method, we see:
\begin{equation}
Q^{-1} \cdot \delta^0_{\epsilon'}  = \begin{bNiceMatrix}
1&  -1 & 1 & 0 & 0 & 0  & 0 \\
0&  1 & 0 & 0 & 0 & 0  & 0 \\
0&  0 & 1 & -1 & 0 & 1  & 0 \\
0&  0 & 0 & 1 & 0 & 0  & 0 \\
0&  0 & 0 & 0 & 1 & 0  & 0 \\
0&  0 & 0 & 0 & 0 & 1  & 0 \\
0&  0 & 0 & 0 & 0 & 0  & 1 \\
\end{bNiceMatrix}  \cdot \begin{bNiceMatrix}
-1&  1 & 0 & 0  & 0 \\
-1&  0 & 1 & 0  & 0\\
0&  -1 & 1 & 0  & 0 \\
0&  -1 & 0 & 1  & 0 \\
0&  -1 & 0 & 0  & 1 \\
0&  0 & -1 & 1  & 0 \\
0&  0 & 0 & -1  & 1 \\
\end{bNiceMatrix} = \begin{bNiceMatrix}
0 & 0 & 0 & 0 & 0\\
-1 & 0 & 1 & 0 & 0 \\
0 & 0 & 0 & 0 & 0 \\
0&  -1 & 0 & 1  & 0 \\
0&  -1 & 0 & 0  & 1 \\
0&  0 & -1 & 1  & 0 \\
0&  0 & 0 & -1  & 1 \\
\end{bNiceMatrix}
\end{equation} 
which is partially reduced, as mentioned in remark \ref{REM: KerIm bases together}. Now we further row-reduce\footnote{Certain care needs to be taken in order for this row-reduction to help us obtain the final result in \eqref{EQ: final result}. In particular, only operations of the form $row_j \mapsto row_j + c \cdot row_k$ for $k < j$.} the matrix until we have established pivot rows:
\begin{align}
\begin{pNiceArray}{c|c}
Q^{-1} \cdot \delta^0_{\epsilon'}& I\\
\end{pNiceArray} = 
\setcounter{MaxMatrixCols}{30}  
&\begin{bNiceArray}{ccccc|ccccccc}
0 & 0 & 0 & 0 & 0 & 1 & 0 & 0 & 0 & 0 & 0& 0 \\
-1 & 0 & 1 & 0 & 0& 0 & 1 & 0 & 0 & 0 & 0& 0 \\
0 & 0 & 0 & 0 & 0 & 0 & 0 & 1 & 0 & 0 & 0& 0 \\
0&  -1 & 0 & 1  & 0& 0 & 0 & 0 & 1 & 0 & 0& 0 \\
0&  -1 & 0 & 0  & 1& 0 & 0 & 0 & 0 & 1 & 0& 0 \\
0&  0 & -1 & 1  & 0& 0 & 0 & 0 & 0 & 0 & 1& 0 \\
0&  0 & 0 & -1  & 1& 0 & 0 & 0 & 0 & 0 & 0& 1 \\
\end{bNiceArray}\\
 \leadsto  &\NiceMatrixOptions{code-for-first-col = \color{blue}}\begin{bNiceArray}[first-col]{ccccc|ccccccc}
v'_1 &0 & 0 & 0 & 0 & 0 & 1 & 0 & 0 & 0 & 0 & 0& 0 \\
v'_2 & -1 & 0 & 1 & 0 & 0& 0 & 1 & 0 & 0 & 0 & 0& 0 \\
v'_3 &0 & 0 & 0 & 0 & 0 & 0 & 0 & 1 & 0 & 0 & 0& 0 \\
v'_4 &0&  -1 & 0 & 1  & 0& 0 & 0 & 0 & 1 & 0 & 0& 0 \\
v'_5 &0&  0 & 0 & -1  & 1& 0 & 0 & 0 & -1 & 1 & 0& 0 \\
v'_6 &0&  0 & -1 & 1  & 0& 0 & 0 & 0 & 0 & 0 & 1& 0 \\
v'_7 &0&  0 & 0 & -1  & 1& 0 & 0 & 0 & 0 & 0 & 0& 1 \\
\CodeAfter
\tikz \draw (2-1) circle (3mm) ;
\tikz \draw (4-2) circle (3mm) ;
\tikz \draw (6-3) circle (3mm) ;
\tikz \draw (7-4) circle (3mm) ;
\end{bNiceArray} = \begin{pNiceArray}{c|c}
P \cdot Q^{-1} \cdot \delta^0_{\epsilon'}& P\\
\end{pNiceArray}
\end{align}
We will come back to how to use the highlighted pivot entries above, as well as the choice of labelings $v_i$ of the rows, in a moment. Our last computational step is to apply the analogous columns operations via $P^{-1}$ (which we solved for by hand by inverting the above row operations).
\begin{align}
 \begin{bNiceMatrix}
\delta^1_{\epsilon'} \cdot Q\\
\hline
Q
\end{bNiceMatrix} \cdot P^{-1}=& \NiceMatrixOptions{code-for-last-col = \color{blue}}\begin{bNiceMatrix}[last-col]
1&  0 & 0 & 0 & 0 & 0  & 0 & u'_{012}\\
0&  0 & 1 & 0 & 0 & 0  & 0 & u'_{123}\\
\hline
1&  1 & -1 & -1 & 0 & 1  & 0& u'_{01}\\
0&  1 & 0 & 0 & 0 & 0  & 0  &u'_{02}\\
0&  0 & 1 & 1 & 0 & -1  & 0&u'_{12}\\
0&  0 & 0 & 1 & 0 & 0  & 0&u'_{13}\\
0&  0 & 0 & 0 & 1 & 0  & 0&u'_{14}\\
0&  0 & 0 & 0 & 0 & 1  & 0&u'_{23}\\
0&  0 & 0 & 0 & 0 & 0  & 1&u'_{34}\\ 
\end{bNiceMatrix} \cdot \begin{bNiceArray}{ccccccc}
1 & 0 & 0 & 0 & 0 & 0 & 0 \\
0 & 1 & 0 & 0 & 0 & 0 & 0 \\
0 & 0 & 1 & 0 & 0 & 0 & 0 \\
0 & 0 & 0 & 1 & 0 & 0 & 0 \\
0& 0 & 0 & 1 & 1 & 0 & 0 \\
0 & 0 & 0 & 0 & 0 & 1 & 0 \\
0 & 0 & 0 & 0 & 0 & 0 & 1 \\
\end{bNiceArray} \\
=& \begin{bNiceMatrix}
\delta^1_{\epsilon'} \cdot Q \cdot P^{-1}\\
\hline
Q \cdot P^{-1}
\end{bNiceMatrix}= \NiceMatrixOptions{code-for-last-col = \color{blue},code-for-last-row = \color{blue}}\begin{bNiceMatrix}[last-col, last-row]
1&  0 & 0 & 0 & 0 & 0  & 0 & u'_{012}\\
0&  0 & 1 & 0 & 0 & 0  & 0 & u'_{123}\\
\hline
1&  1 & -1 & -1 & 0 & 1  & 0& u'_{01}\\
0&  1 & 0 & 0 & 0 & 0  & 0  &u'_{02}\\
0&  0 & 1 & 1 & 0 & -1  & 0&u'_{12}\\
0&  0 & 0 & 1 & 0 & 0  & 0&u'_{13}\\
0&  0 & 0 & 1 & 1 & 0  & 0&u'_{14}\\
0&  0 & 0 & 0 & 0 & 1  & 0&u'_{23}\\
0&  0 & 0 & 0 & 0 & 0  & 1&u'_{34}\\ 
v'_1&  v'_2 & v'_3 & v'_4 & v'_5 & v'_6  & v'_7
\CodeAfter
\tikz \draw[rounded corners] (3-2) -| (last-|3) -- (last-|2) |- (3-2) ;
\tikz \draw[rounded corners] (3-4) -| (last-|5) -- (last-|4) |- (3-4) ;
\tikz \draw[rounded corners] (3-5) -| (last-|6) -- (last-|5) |- (3-5) ;
\tikz \draw[rounded corners] (3-6) -| (last-|7) -- (last-|6) |- (3-6) ;
\tikz \draw[rounded corners] (3-7) -| (last-|8) -- (last-|7) |- (3-7) ;
\end{bNiceMatrix}
\end{align}
Using the same analysis as \eqref{EQ:BasisKerDelta0prime}, the columns of $Q \cdot P^{-1}$ directly below the zero columns of $\delta^1_{\epsilon'} \cdot Q \cdot P^{-1}$ provide a basis for the null space:
\begin{align}\label{EQ:BasisKerDelta1prime}
ker\left( \delta^1_{\epsilon'} \right) = &\text{span} \left\{v'_2= u'_{01}+ u'_{02}, v_4 = -u'_{01} + u'_{12} + u'_{13} + u'_{14}, v'_5 = u'_{14},\right.\\
&\left. v'_6 = u'_{01}- u'_{12} + u'_{23} , v'_7 =  u'_{34}\right\}. \nonumber
\end{align}
Now the pivots from $P \cdot Q^{-1} \delta^0_{\epsilon'}$ will tell us which subset of this null space basis set forms a basis for the image subpace:
\begin{equation}\label{EQ:BasisImDelta0prime}
im\left( \delta^0_{\epsilon'} \right) = \text{span} \left\{v'_2, v'_4, v'_6, v'_7   \right\}. 
\end{equation}

Recalling that if $U \subset W$ are vector spaces with respective bases $\{u_1, \dots, u_n\}$ and $\{u_1, \dots, u_n, w_1, \dots w_k\}$, then we have a basis for the quotient space $W/U$ given by $\{ [w_1], \dots, [w_k] \}$. Using this fact along with the bases from \eqref{EQ:BasisKerDelta0prime}, \eqref{EQ:BasisKerDelta1prime}, and \eqref{EQ:BasisImDelta0prime} we finally obtain the bases for the vector spaces:

\begin{align}
\check{H}^1_{\epsilon'} &= ker\left( \delta^1_{\epsilon'} \right) / im\left( \delta^{0}_{\epsilon'}\right) \cong  \text{span} \left\{[v'_5] = [u'_{14}] \right\}\cong \mathbb{R}^1 \label{EQ:BasesCohomPrime}\\
\check{H}^0_{\epsilon'} &= ker\left( \delta^0_{\epsilon'} \right) / im\left( \delta^{-1}_{\epsilon'}\right) \cong  \text{span} \left\{u'_0 + u'_1 + u'_2 + u'_3 + u'_4 \right\} / \{0\} \nonumber \\
&\cong  \text{span} \left\{u'_0 + u'_1 + u'_2 + u'_3 + u'_4 \right\} \cong \mathbb{R}^1 \nonumber
\end{align} 
The dimensions of $\check{H}^0_{\epsilon'}$ and $\check{H}^1_{\epsilon'}$ correspond to the Betti numbers we would expect for $\mathcal{D}_{\epsilon'}$: there is a single connected component so we have $\beta^0 = dim \left( \check{H}^0_{\epsilon'} \right) = 1$ and there is a single non-trivial loop so we have $\beta^1 = dim \left( \check{H}^1_{\epsilon'} \right) = 1$. 

\subsection{A basis for the cohomology induced by $\epsilon$, compatible with that of $\epsilon'$.}

Next, we want to compute bases for the column and null spaces of the boundary maps $\delta^0_{\epsilon = 1.4}, \delta^1_{\epsilon = 1.4},$ etc. Once again, we must be careful to compute these bases in such a way that they are coherent with the previous bases. This means that we have to ensure the basis for each $im\left( \delta^{p}_{\epsilon}\right)$ includes into the basis we use for each $ker\left( \delta^{p+1}_{\epsilon}\right)$ just like above. However, there is an additional level of difficulty: we must also keep in mind the bases used for computing the cohomology with respect to $\epsilon'$ so that we can discuss which classes ``persist''. If we are careful (i.e. using an algorithm we are trying to share via example) we can proceed with our reductions in a way that all of the above is satisfied. 

The essential idea is to imagine that $\epsilon'$ is the largest of the finite list of resolutions which produce different complexes. We maintain all of the components, but fill in zeros for basis elements which vanish under the new $\epsilon$. This allows us to use the same row/column  transformations and simply read off the resulting information. To start we apply our column transformation $L$ from \eqref{EQ: H0 epsilon prime reduction} to obtain a basis for $ker(\delta^0_{\epsilon})$:

\begin{align}
\begin{bNiceMatrix}
\delta^0_{\epsilon}\\
\hline
Id_{\check{C}^0_{\epsilon}}
\end{bNiceMatrix}  \cdot L= &
\NiceMatrixOptions{code-for-first-row = \color{blue},
code-for-last-col = \color{blue}}\begin{bNiceMatrix}[first-row]
u_0 & u_1 & u_2 & u_3 & u_4\\
-1&  1 & 0 & 0  & 0 \\
-1&  0 & 1 & 0  & 0 \\
0&  -1 & 1 & 0  & 0 \\
0&  -1 & 0 & 1  & 0 \\
0&  0 & 0 & 0  & 0 \\
0&  0 & -1 & 1  & 0 \\
0&  0 & 0 & 0  & 0 \\
\hline
1&  0 & 0 & 0  & 0 \\
0&  1 & 0 & 0  & 0 \\
0&  0 & 1 & 0  & 0 \\
0&  0 & 0 & 1  & 0 \\
0&  0 & 0 & 0  & 1
\end{bNiceMatrix} \cdot \begin{bNiceMatrix}
1&  1 & 1 & 1  & 1 \\
0&  1 & 1 & 1  & 1\\
0&  0 & 1 & 1  & 1\\
0&  0 & 0 & 1  & 1 \\
0&  0 & 0 & 0  & 1 
\end{bNiceMatrix}=\NiceMatrixOptions{code-for-first-row = \color{blue},
code-for-last-col = \color{blue}}\begin{bNiceMatrix}[last-col]
-1&  0 & 0 & 0  & 0 & u_{01}\\
-1&  -1 & 0 & 0  & 0& u_{02}\\
0&  -1 & 0 & 0  & 0& u_{12}\\
0&  -1 & -1 & 0  & 0& u_{13}\\
0&  0 & 0 & 0  & 0& u_{14}\\
0&  0 & -1 & 0  & 0  &u_{23}\\
0&  0 & 0 & 0  & 0& u_{34}\\
\hline
1&  1 & 1 & 1  & 1 & u_{0}\\
0&  1 & 1 & 1  & 1& u_{1}\\
0&  0 & 1 & 1  & 1& u_{2}\\
0&  0 & 0 & 1  & 1 & u_{3}\\
0&  0 & 0 & 0  & 1 & u_{4}
\CodeAfter
\tikz \draw[rounded corners] (8-5) -| (last-|6) -- (last-|5) |- (8-5) ;
\tikz \draw[rounded corners] (8-4) -| (last-|5) -- (last-|4) |- (8-4) ;
\end{bNiceMatrix} = \begin{bNiceMatrix}
\delta^0_{\epsilon} \cdot L\\
\hline
Id_{\check{C}^0_{\epsilon}} \cdot L
\end{bNiceMatrix}  
\end{align}
Once again, the circled columns provide a basis for the kernel of $\delta^0_{\epsilon}$: 
\begin{align}\label{EQ:BasisKerDelta0}
ker\left( \delta^0_{\epsilon} \right) &= \text{span} \left\{u_0 + u_1 + u_2 + u_3 , u_0 + u_1 + u_2 + u_3 + u_4 \right\} . 
\end{align}
Next, we apply our $P$ and $Q^{-1}$ row-transformations to $\delta^0_{\epsilon}$ from \eqref{EQ:delta0example},
\begin{align}
P \cdot Q^{-1} \cdot \delta^0_{\epsilon} &= \NiceMatrixOptions{code-for-first-col = \color{blue}}\begin{bNiceArray}[first-col]{ccccc}
v_1 & 0 & 0 & 0 & 0 & 0  \\
v_2 & -1 & 0 & 1 & 0 & 0  \\
v_3 & 0 & 0 & 0 & 0 & 0 \\
v_4 & 0 & -1 & 0 & 1 & 0 \\
v_5 & 0 & 1 & 0 & -1 & 0  \\
v_6 & 0 & 0 & -1 & 1 & 0  \\
v_7 & 0 & 0 & 0 & 0 & 0\\
\end{bNiceArray}
\end{align} 
en route to obtain a basis for the image. We will come back to this matrix shortly. There are also the corresponding column-transformations on $\delta^1_{\epsilon}$ from \eqref{EQ:delta1example},
\begin{align}
 \begin{bNiceMatrix}
\delta^1_{\epsilon'} \\
\hline
I_{\epsilon}
\end{bNiceMatrix} \cdot (Q \cdot P^{-1}) =& \NiceMatrixOptions{code-for-last-col = \color{blue},code-for-last-row = \color{blue}}\begin{bNiceMatrix}[last-col, last-row]
1&  -1 & 1 & 0 & 0 & 0  & 0 & u_{012}\\
0&  0 & 0 & 0 & 0 & 0  & 0 & u_{123}\\
\hline
1&  0 & 0 & 0 & 0 & 0  & 0 & u_{01}\\
0&  1 & 0 & 0 & 0 & 0  & 0 & u_{02}\\
0&  0 & 1 & 0 & 0 & 0  & 0 & u_{12}\\
0&  0 & 0 & 1 & 0 & 0  & 0 & u_{13}\\
0&  0 & 0 & 0 & 0 & 0  & 0 & u_{14}\\
0&  0 & 0 & 0 & 0 & 1  & 0 & u_{23}\\
0&  0 & 0 & 0 & 0 & 0  & 0 & u_{34}\\  
v_1&  v_2 & v_3 & v_4 & v_5 & v_6  & v_7
\CodeAfter
\end{bNiceMatrix}\\
 =  &\begin{bNiceMatrix}
\delta^1_{\epsilon}  \cdot I_{\epsilon} \cdot Q \cdot P^{-1}\\
\hline
 I_{\epsilon}  \cdot Q \cdot P^{-1}
\end{bNiceMatrix}= \NiceMatrixOptions{code-for-last-col = \color{blue},code-for-last-row = \color{blue}}\begin{bNiceMatrix}[last-col, last-row]
1&  0 & 0 & 0 & 0 & 0  & 0 & u_{012}\\
0&  0 & 0 & 0 & 0 & 0  & 0 & u_{123}\\
\hline
1&  1 & -1 & -1 & 0 & 1  & 0 & u_{01}\\
0&  1 & 0 & 0 & 0 & 0  & 0 &u_{02}\\
0&  0 & 1 & 1 & 0 & -1  & 0&u_{12}\\
0&  0 & 0 & 1 & 0 & 0  & 0&u_{13}\\
0&  0 & 0 & 0 & 0 & 0  & 0&u_{14}\\
0&  0 & 0 & 0 & 0 & 1  & 0&u_{23}\\
0&  0 & 0 & 0 & 0 & 0  & 0&u_{34}\\  
v_1&  v_2 & v_3 & v_4 & v_5 & v_6  & v_7
\end{bNiceMatrix}
\end{align}
In the above matrix, the columns in $ I_{\epsilon}  \cdot Q \cdot P^{-1}$ beneath the trivial columns of $\delta^1_{\epsilon}  \cdot I_{\epsilon} \cdot Q \cdot P^{-1}$ provide us the spanning set $\{v_2, v_3, v_4,v_6\}$ for the kernel. The point of doing things this way is that our chain map sends our previous basis for $ker(\delta^1_{\epsilon'})$ to this spanning set for $ker(\delta^1_{\epsilon})$. However, these vectors are not necessarily linearly independent, so we could first check the pivots of the submatrix formed by these columns, to find a basis for this spanning space, and then find the pivots of the corresponding submatrix of the transformed $\delta^0$ above. In our case, our vectors are linearly independent and so, 
\begin{align}
ker\left( \delta^1_{\epsilon} \right) = &\text{span} \left\{v_2 = u_{01} + u_{02}, v_3 = -u_{01} + u_{12}, v_4 = -u_{01} + u_{12} + u_{13}, \right. \nonumber\\
&\left.  v_6 = u_{01} - u_{12} + u_{23}  \right\}. 
\end{align}
To see which of these vectors form a basis in the image (note $f$ of the image spans the image!), we can look at the pivots of $P Q^{-1} \delta^0$:
\begin{align}
P \cdot Q^{-1} \cdot \delta^0_{\epsilon} &= \NiceMatrixOptions{code-for-first-col = \color{blue}}\begin{bNiceArray}[first-col]{ccccc}
v_1 & 0 & 0 & 0 & 0 & 0  \\
v_2 & -1 & 0 & 1 & 0 & 0  \\
v_3 & 0 & 0 & 0 & 0 & 0 \\
v_4 & 0 & -1 & 0 & 1 & 0 \\
v_5 & 0 & 1 & 0 & -1 & 0  \\
v_6 & 0 & 0 & -1 & 1 & 0  \\
v_7 & 0 & 0 & 0 & 0 & 0\\
\end{bNiceArray} \leadsto & \begin{bNiceArray}{c} v_2 \\ v_3 \\ v_4 \\ v_6  \end{bNiceArray} = &\NiceMatrixOptions{code-for-first-col = \color{blue}}\begin{bNiceMatrix}[first-col]
v_2 & -1 & 0 & 1 & 0 & 0  \\
v_3 & 0 & 0 & 0 & 0 & 0 \\
v_4 & 0 & -1 & 0 & 1 & 0 \\
v_6 & 0 & 0 & -1 & 1 & 0  \\
\end{bNiceMatrix}\\
\intertext{and now applying column reduction (to see which rows are independent!),}
 \xmapsto{\substack{c_3: c_3+ c_1\\ c_4: c_4 + c_2 \\c_4: c_4 + c_3 }}&\begin{bNiceMatrix} 
 -1 & 0 & 0 & 0 & 0  \\
 0 & 0 & 0 & 0 & 0 \\
 0 & -1 & 0 & 0 & 0 \\
 0 & 0 & -1 & 0 & 0  \\
\CodeAfter
\tikz \draw (1-1) circle (3mm) ;
\tikz \draw (3-2) circle (3mm) ;
\tikz \draw (4-3) circle (3mm) ;
\end{bNiceMatrix}
\end{align}
and so
\[ im\left( \delta^0_{\epsilon'} \right) = \text{span} \left\{v_2, v_4, v_6   \right\}.\]
In summary we have:
\begin{align}
\check{H}^0_{\epsilon} &= ker\left( \delta^0_{\epsilon} \right) / im\left( \delta^{-1}_{\epsilon}\right) \cong  \text{span} \left\{u_0 + u_1 + u_2 + u_3, u_0 + u_1 + u_2 + u_3 + u_4  \right\} / \{0\} \nonumber \\
&\cong  \text{span} \left\{u_0 + u_1 + u_2 + u_3, u_0 + u_1 + u_2 + u_3 + u_4 \right\} \cong \mathbb{R}^2 \nonumber\\
\check{H}^1_{\epsilon} &= ker\left( \delta^1_{\epsilon} \right) / im\left( \delta^{0}_{\epsilon}\right) \cong  \text{span} \left\{[v_3] = [-u_{01} + u_{12}] \right\}\cong \mathbb{R}^1
\end{align} 
The dimensions of the two above spaces $\check{H}^0_{\epsilon}$ and $\check{H}^1_{\epsilon}$ correspond to the Betti numbers we would expect for $\mathcal{D}_{\epsilon}$: there are two connected components so we have $\beta^0 = dim \left( \check{H}^0_{\epsilon} \right) = 2$ and there is a single non-trivial loop so we have $\beta^1 = dim \left( \check{H}^1_{\epsilon} \right) = 1$. 

We present a summary of all of these calculations so far in the following diagram:
\begin{equation}\label{EQ: final result}
\begin{tikzpicture}
 \foreach \X/\z in {0/-1,1/0,2/1}{
	\foreach \Y in {0,1}{
 \draw (5*\X, 5*\Y) rectangle (5*\X+4, 5*\Y+3);
 \draw[->] (5*\X ,5*\Y) rectangle (5*\X+3.8 ,5*\Y+ 2.5);
 \node[below left] at (5*\X+3.8 ,5*\Y+ 2.5) {\tiny{$ker\left( \delta^{\X} \right)$}};

 }}
 \foreach \X/\z in {1/0,2/1}{
	\foreach \Y in {0,1}{
 \draw[->] (5*\X ,5*\Y) rectangle (5*\X+3.5 ,5*\Y+ 1.5);
  \node[below left] at (5*\X+3.5 ,5*\Y+ 1.5) {\tiny{$im\left( \delta^{\z} \right)$}};
 }}
  \foreach \X in {0,1,2,3}{
 	\foreach \Y in {0, 1}{
	\draw[->] (5*\X -0.75,5*\Y+ 1.5) to (5*\X -0.25,5*\Y+ 1.5);
	}}
\foreach \X in {0,1,2}{
 \draw[->] (5*\X +1.5,4.5) to (5*\X +1.5,3.5)node [label=right:{$f^{\X}$}, yshift=0.23cm, xshift=-0.1cm] {};}
 \node[above] at (4.5, 6.5) {$\delta^0_{\epsilon'}$};
 \node[above] at (9.5, 6.5) {$\delta^1_{\epsilon'}$};
 \node[above] at (4.5, 1.5) {$\delta^0_{\epsilon}$};
 \node[above] at (9.5, 1.5) {$\delta^1_{\epsilon}$};
\foreach \X in {-1,15}{
	\foreach\Y in {1.5,6.5}{
	\node at (\X, \Y) {$0$};}}
\foreach \X in {0,1,2}{
	\node at (5*\X+0.3, 7.75) {$\check{C}^{\X}_{\epsilon'}$};}
\foreach \X in {0,1,2}{
	\node at (5*\X+0.3, 2.75) {$\check{C}^{\X}_{\epsilon}$};}
\node at (1.85,5.85) {\tiny{$u'_0 +  u'_1 + u'_2 + u'_3 + u'_4 $}};
\node at (1.85,0.85) {\tiny{$u_0  + u_1 + u_2 + u_3 $}};
\node at (1.85,0.55) {\tiny{$u_0  + u_1 + u_2 + u_3+u_4 $}};
\node at (5.85,6.85) {\tiny{$u'_{14}$}};
\node at (6.2,6.2) {\tiny{$u'_{01}+ u'_{02},$}};
\node at (6.7,5.9) {\tiny{$-u'_{01} + u'_{12} +u'_{13} + u'_{14},$}};
\node at (6.5,5.6) {\tiny{$u'_{01} - u'_{12} + u'_{23}$}};
\node at (6.7,5.3) {\tiny{$u'_{34}$}};
\node at (6.2,1.85) {\tiny{$-u_{01} + u_{12}$}};
\node at (6.75,0.85) {\tiny{$u_{01} + u_{02}$}};
\node at (6.2,0.55) {\tiny{$-u_{01}+u_{12}+u_{13}$}};
\node at (6.2,0.25) {\tiny{$u_{01}-u_{12}+u_{23}$}};
\node at (11.2,6.85) {\tiny{$0$}};
\node at (11.2,5.85) {\tiny{$u'_{012}, u'_{123}$}};
\node at (11.2,1.85) {\tiny{$0$}};
\node at (11.2,0.85) {\tiny{$u_{012}$}};
\end{tikzpicture}
\end{equation}
\subsection{Explicit persistence from $\epsilon'$ to $\epsilon$.}
The final goal is to use these calculations to study the \emph{persistence} of these topological features (clusters and loops) as the radii of our discs $\epsilon \le \epsilon'$ vary. Let's start with $\check{H}^0_{\epsilon'} \xrightarrow{f} \check{H}^0_{\epsilon}$ . Intuitively, we see from Figure \ref{fig: Simple Point Cloud} that the two separate connected components for $\epsilon$ \emph{merge} together as the radius increases to $\epsilon'$. Algebraically, this is seen by: (i) the class $c'=[u'_0 + \dots + u'_4] \in \check{H}^0_{\epsilon'}$ gets sent by $f^0$ to the class $c=[u_0 + u_1 + u_2 + u_3,u_4] \in \check{H}^0_{\epsilon}$, and so the class $c$ ``(co)-persists'' while (ii) the class $[u_0 + u_1 + u_2 + u_3] \in \check{H}^0_{\epsilon}$ ''(co)-dies'' as it is not in the image of a class from $\epsilon'$. In practice some more care is taken to decide whom merges which whom if there are multiple (co)-deaths.

Next for $\check{H}^1_{\epsilon'} \xrightarrow{f} \check{H}^1_{\epsilon}$, the cohomology computations tell us that each point cloud has a non-trivial loop at the chosen radius: $[u'_{14}] \in \check{H}^1_{\epsilon'}$ and $[-u_{01} +u_{12}] \in \check{H}^1_{\epsilon}$  . However, we can see from Figure \ref{fig: Simple Point Cloud} that these two loops should not be considered the same. In persistent (co)homology, we would say just by looking at the figure that the loop-class $[u'_{14}]$ is ``(co)-born'' at $\epsilon'$ while loop-class $[-u_{01}+u_{12}]$ ``(co)-dies'' at $\epsilon'$. Algebraically, this is seen by $f^1\left( [u'_{14}]\right) = 0 \in \check{H}^1_{\epsilon}$ while $[-u_{01} + u_{12}]$ is not in the image of $f^1$.

\end{document}